\documentclass{IEEEtran4PSCC}
\ifCLASSINFOpdf
   \usepackage[pdftex]{graphicx}
   \graphicspath{{./Images/}}
   \DeclareGraphicsExtensions{.pdf,.jpeg,.png}  
\else
\fi

\usepackage[cmex10]{amsmath}
\usepackage{algorithmic}
\usepackage{array}
\usepackage{url}
\usepackage{hyperref}
\usepackage{graphicx}
\usepackage{color}

\AtBeginDocument{
  \DeclareSymbolFont{AMSb}{U}{msb}{m}{n}
  \DeclareSymbolFontAlphabet{\mathbb}{AMSb}}

\begin{document}
\title{Modeling analysis and optimization for European network data merging}

\author{
\IEEEauthorblockN{
	Manuel Ruiz\\
	Othman Moumni Abdou\\
	Arnaud Renaud
}
\IEEEauthorblockA{ARTELYS\\
Paris, France\\
firstname.lastname@artelys.com}
\and
\IEEEauthorblockN{
	Jean Maeght\\
	Mireille Lefevre\\
	Patrick Panciatici
}
\IEEEauthorblockA{RTE\\
Versailles, France\\
firstname.lastname@rte-france.com}
}

\maketitle
\begin{abstract}
In this paper, the problem of building a consistent European network state based on the data provided by different Transmission System Operators (TSOs) is addressed. 
A hierarchical merging procedure is introduced and consists in the resolution of several Optimal Power Flow problems (OPFs).
Results on the European network demonstrate the interest of this procedure on real-life cases and highlight the benefits of using a hierarchical multi-objective approach.

\end{abstract}


\section{Introduction}
\label{sec:intro}
One of the main concern of European Transmission System Operators in Europe is keeping the system in a secure state. Security assessment tools help to assess risk in the system and provide preventive and curative actions when necessary. They are performed ex post and require a description of the network at an earlier time $t$.
In an ideal case, each TSO would provide \emph{accurate} data, representing the network state at the \emph{exact time}~$t$.
A simple power flow computation would then be sufficient to establish the state of the whole network.

Unfortunately, available data are not synchronized between TSOs, input data are snapshot generated each 15 minutes that are always slightly non-synchronized. They can also be missing; in this case they are replaced by the previous snapshot or the forecast made the day before to have a complete description of the network.
Thus, non-synchronized and erroneous data have to be collected in real-time and merged in order to build a consistent state of the European network for security assessment. Additionally, the merging must limit the impact of the erroneous data on the consolidation of reliable snapshot data.

The problem of state estimation enhancement~\cite{dancre2002optimal} leads to the resolution of Optimal Power Flow where the objective is to find a feasible network state minimizing the deviations from the well-estimated parts of the network.
This problem has already been considered by RTE, the French TSO, who introduced an objective of minimizing a weighted
least squares distance and implemented the OPF in an optimizer for merging European networks.
This merging tool is used at {\color{blue} \href{http://www.coreso.eu/}{CORESO}}
(COoRdination of Electricity System Operators) on a day-to-day basis since 2009. More details on power system optimization and general OPF formulations can be found in~\cite{Cain12},~\cite{Fliscounakis2006}. The work presented in this paper represents possible way of improvement of the CORESO merging tool with the aim to use a multi-objective hierarchical approach. A first study~\cite{pscc2015} have already been published in this way to evaluate the impact of the deviation measure in the optimization process.

\vskip.3cm
The paper is organized as follows. Section~\ref{sec:context} presents the context of the study and why the {\color{blue} \href{http://www.itesla-project.eu/}{iTesla}}  project leads to this problematic. After introducing the mathematical formulation of the modified ACOPF in section~\ref{sec:math}, main contributions are detailed in section~\ref{sec:main_contrib}. The hierarchical merging methodology is presented in section~\ref{sec:hierachical_approach}. Finally, computational experiments on large scale power system network are presented in section~\ref{sec:exp}.

\section{Context}
\label{sec:context}

\subsection{The iTesla project}
\label{subsec:itesla_project}

The Innovative Tools for Electrical System Security within Large Areas (\color{blue} \href{http://www.itesla-project.eu/}{iTesla}\color{black}) project is
a European research and development project co-funded by the European Commission~\cite{iTeslaWS}
and involving 6 European TSOs, 13 other R\&D providers and several academic partners.
The aim of this project is to develop a flexible inter-operable toolbox which is able to support
future operations of the pan-European electricity transmission network.
This toolbox shall enable to perform:
\begin{itemize}
	\item Accurate security assessments taking system dynamics into account using time-domain simulations.
	\item Risk assessments taking into account: uncertainties (especially intermittent power generation uncertainties),
	probabilities of contingencies and possible failures of corrective actions.
\end{itemize}

\vskip.3cm
The iTesla toolbox includes two main components:
\begin{itemize}
	\item an off-line platform, running every week, in charge of exploring the network
	state space to draw the separation between stable and unstable states (using data
	mining techniques and high-fidelity dynamic simulations). This platform is especially in
	charge of identifying security rules to apply.
	\item an on-line platform, running every day, dedicated to security assessment.
	This second platform evaluates the off-line platform security rules on the current
	network and issues curative or preventive actions, when necessary.
\end{itemize}

\begin{figure}[!h]
\centering
\includegraphics[width=3.5in]{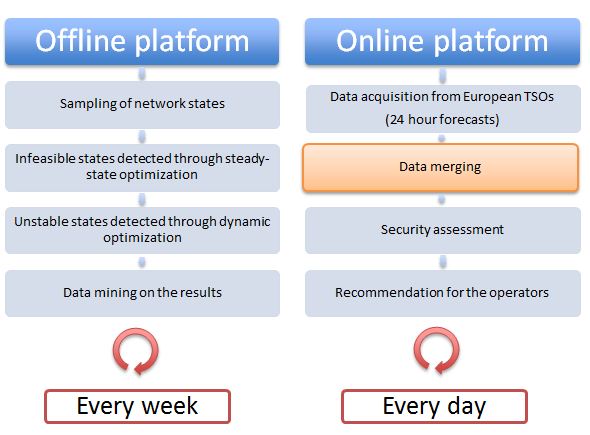}
\caption{The iTesla workflow}
\label{itesla_workflow}
\end{figure}

The whole iTesla workflow is described in the figure~\ref{itesla_workflow}. The on-line platform is also in charge of acquiring and merging data provided by European TSOs. The work presented in this paper is part of data merging module presented in the following subsection.

\subsection{Data merging}
\label{subsec:data_merging}
The data merging module is part of the online platform which aims to consolidate TSOs' networks into a single, consistent, merged European network. Result of this module is used as a description of the underlying network to perform security assessments and provide recommendation for TSOs to keep the system in a secure state.

Input data are the result of a topological merge of European TSOs' network data performed with the most reliable data available fot time $t$. The output is a consistent European network state for the time $t$, minimizing the deviation from the data provided by TSOs.

In the data, the network state for each country is consistent but the whole network is not.
One of the reasons is the power flow on interconnections:
when a European AC load flow is solved, the active power flow on interconnections is computed using the characteristics of the interconnections and the voltage level of the buses. Computed and measured values are usually different since data from the different TSOs are not perfectly synchronized.

The data provided by each TSO become targets and the aim is to find a feasible merged network state minimizing the deviations from target data. Deviations that are minimized are related to voltage level of bus, active power flow on interconnections and loads.

Additionally, there are different types of data and different levels of confidence:
\begin{itemize}
 \item Time differences: the aim is to compute a network state at time $t$ but the data provided by the TSOs
 can be up to one hour older
 \item There are different reliability levels:
 \begin{enumerate}
  \item Snapshots (SN): data based on real-time observations and measurements on the network.
    If they are recent, these data are very reliable.
  \item Day Ahead Congestion Forecasts (DACF): these are forecasts data for a given time.
    If a TSO is not able to provide its snapshot data in time, the DACF data (which is
    provided the day before) are used as a replacement.
    Obviously these data are less reliable than SN data. Although, using this data must not impact the merging of the SN part of the network; erroneous data must not be spread to the other parts of the network.
 \end{enumerate}
\end{itemize}

\vskip.3cm
Consequently, the merging procedure involves multiple minimization objectives depending on deviation type (voltage levels, active power flows or loads) and data reliability level (SN or DACF). In the following section, the mathematical formulation of network transmission data consolidation is presented.

\section{Mathematical formulation}
\label{sec:math}
The power flow is modeled using an AC formulation. The model is similar to the polar PQV ACOPF presented in ~\cite{Cain12}.
Additionally, the loads are decomposed into two parts: a constant part given by TSOs plus a slack variable for the deviation. 
The deviations model the gap between the values provided by the TSOs and the values in a feasible solution.
Two types of buses are distinguished: PV and PQ.
For PV buses, active production, voltage and load are given. The active production is fixed and the voltage level and load are used as targets.
For PQ buses, active and reactive production and load are specified. The active and reactive production are fixed and load is used as target data.
The voltage of PQ buses and the reactive power of PV buses are within their bounds but do not appear in the objectives of optimization models.

In order to measure the deviations from the collected data, different penalty functions can be considered.
The study in~\cite{pscc2015} has evaluated the influence of the norms used as penalty functions in network merging procedure. The results of the study showed that the deviation measures can have a significant impact on how the data are consolidated and how the deviations are spread over the buses and the interconnections. 

The deviations are measured using the Manhattan norm, the sum of squares and Huber function. Choice of these penalty functions is based on the results of the study in~\cite{pscc2015}. 

\vskip.3cm
The Manhattan norm ($L_1$) sums the absolute values of the deviations. It helps focusing on minimizing the total sum of the deviations.
Additionally, since the impact is linear there is no advantage at modifying several variables preferably to a single variable.

This norm is used to measure the deviation of loads for SN buses ($obj_{load}^{sn}$). As SN data are assumed to have a good confidence level, it ensures that deviations of loads are focused on a small number of buses.
Moreover, it helps to localize problems and analyze the merged network state.

\vskip.3cm
The sum of squares ($L_2$) favors small modifications over many variables than significant modification on few ones.
This prevents the deviations from being too wide.
It helps finding a feasible solution within a limited range of the original one which is especially interesting
if TSOs data are reliable but can fail to correct a localized error.

This norm is used to measure deviation of voltage levels for PV buses ($obj_{V}$) as these data are very reliable. It helps to find feasible solution within a limited range of the original voltage levels.

It is also used as penalty function for deviation of loads for DACF buses ($obj_{load}^{dacf}$). Since DACF data are used as replacement of SN data and are less reliable, all loads from DACF data are potentially erroneous. Using the sum of squares as penalty function allows to spread deviations over DACF buses.

\vskip.3cm
The Huber function is quadratic for small differences and linear for wide variations.
This favors spreading deviations provided they remain small and does not strongly penalize a significant change of a variable. It combines both the advantages of $L_1$ and $L_2$ norms. This approach has already been successfully applied in~\cite{Fliscounakis2006}. Additionally, the sum of two Huber functions is constant and minimum over a confidence interval. This is especially interesting if there are different targets for the same data as it is the case for active power flow on interconnections where each neighboring TSO provides a target for a given interconnection. 

Therefore, this function is used as penalty function for the deviation of active flows on interconnections ($obj_{interco}$) as it helps take into account two different targets. The variation of active power flow between these two targets have no impact on $obj_{interco}$ unless the variation is wider.

\section{Main contribution}
\label{sec:main_contrib}
The four deviation minimization objectives presented in section~\ref{sec:math} may be conflicting. For example, minimizing the deviation of load for a SN bus connected to an interconnection could increase the deviation of active power flow on this interconnection. In practice, there does not exist a single solution that simultaneously optimizes each objective and thus a trade-off exists between the different objectives.

Classical merging methods, like the one used at {\color{blue} \href{http://www.coreso.eu/}{CORESO}} uses a single-objective function problem that involves, separately, all the objectives to minimize deviations.
This method normalizes the objectives into one single-objective by multiplying each objective with a weight. The set of weights used in such objective could be based on a prioritization of all types of deviations but it will highly influences the quality of the merged network state obtained by spreading deviations over the objectives.

Hence, one can consider using a hierarchical approach by solving successive ACOPFs using different objectives at each step. This method seeks the appropriate trade-off between the different objectives in order to find a feasible network state using the order of preference of the objectives.
Moreover, with the advent of more efficient non-linear solvers, considering a hierarchical merging approach has become possible as it is fast enough to be used every 15 minutes in an operational context.

The methodology presented in this paper is successfully applied to the very high voltage European transportation network. Results of tests on real operational data are presented in section~\ref{sec:exp}. In the subsequent section~\ref{sec:hierachical_approach}, the hierarchical approach is detailed.

\section{Hierarchical approach}
\label{sec:hierachical_approach}
The hierarchical merging procedure consists of solving several successive modified ACOPFs.

One of the main difficulties of this method is how to constraint the problem at each step to the locally optimal solution space of the preferred objectives (objectives of previous steps). In fact, all objectives have to be minimized in a manner such that they do not cause an increase in any of the more preferred objectives. This rule has been relaxed to slightly expand the solution space in order to avoid numerical difficulties as OPF are highly non-convex continuous optimization problems.

\vskip.3cm
In the first step, only the $obj_{interco}$ is used as objective function for the ACOPF model. 
Since only deviations of active power flows on interconnections are minimized, other data (voltage levels and loads) may be modified without penalization, while respecting the constraints of the ACOPF model. The resulting flows are thus generally within the two targets provided by the TSOs.
After this first step, active power flows on interconnections are bounded in the interval containing the two targets and the value obtained in this first step. A tolerance is used to extend the range of accepted values and accept neighboring values.

The second step aims to minimize $obj_{V}$, starting from network state obtained in the first step. The resulting merged network state has deviations mainly concentrated on loads since active power flows on interconnections are bounded and deviations on voltage levels of PV buses are minimized.
After this second step, voltage levels on PV buses are fixed using a small tolerance. 

In the third step, $obj_{load}^{sn}$ is used to minimize the deviations of loads from SN data. The optimization tends to prior modification of DACF loads over SN loads, since there is no penalization associated to loads from DACF data. In the SN perimeter, the use of the $L_1$ norm helps to localize deviations and avoid making several small deviations.
At the end of this step, deviations of loads from SN data are fixed.

In the last step, the only degrees of freedom in the model are loads from DACF data. Flows on interconnections and voltage levels are fixed using tolerances while loads in SN perimeter are fixed. This last step focuses on spreading the deviations of loads on DACF buses by  minimizing the objective $obj_{load}^{dacf}$.

At the end of the hierarchical merging procedure, a consistent merged network state is obtained where deviations from collected data are minimized using the adequate trade-off between the different objectives. The deviations are prioritized using the following order: DACF loads (lowest priority), SN loads, voltage levels of PV buses and then active power flow on interconnections (highest priority). 

\section{Computational experiments}

\label{sec:exp}
The iTesla toolbox will be operating at a Pan-European level. In this work, data are collected from 15 TSOs in 11 countries. The grid is composed of around 6000 buses, 8500 branches, 3500 loads, 900 production units and about 180 interconnections between the different TSOs' networks. This underlines how challenging the project is.

The hierarchical merging procedure is implemented using AMPL~\cite{citeulike:761822} and solved with KNITRO~\cite{Byrd99aninterior, Byrd2000}.

\subsection{Merging indicators}
In order to evaluate the quality of the resulting merged network state, the following indicators are defined to estimate the deviation from TSOs data:
\begin{itemize}
\item \textit{Interco SN\_SN}: Number of interconnections between two TSOs with SN data where the active power flows are not within the two targets provided by two neighboring TSOs.
\item \textit{Interco SN\_DACF}: Number of interconnections between a TSO with SN data and a TSO with DACF data such that the active power flows target in the SN data is not satisfied.
\item \textit{V SN}: Number of SN buses where voltage level targets are not met.
\item \textit{V DACF}: Number of DACF buses where voltage level targets are not satisfied.
\item \textit{Loads ACT SN}: Number of SN buses where active load targets are not met.
\item \textit{Loads REA SN}: Number of SN buses where reactive load targets are not satisfied.
\end{itemize}
Acceptable deviation tolerances are used in order to only consider real deviations. The deviation of active power flow on interconnection is non-zero if its absolute value is greater than $5MW$. Concerning voltage levels, only absolute deviations greater than $0.1kV$ are taken into account. Deviations of active loads and reactive loads are non-zero if the absolute value is respectively greater than $1MW$ and $1Mvar$.

\subsection{Experimentation to determine the modeling approach}
\label{subsec:first_exp}

First experiments allow us to validate the hierarchical merging procedure and ensure that the method gives good merging results on large scale networks. Several modeling approach were tested during these first experiments which allows us to determine the best modeling alternative: choice of penalty functions, definition of the objectives and their hierarchical order, numerical tolerances to avoid infeasibilities, etc. Moreover, a load flow is used in the iTesla toolbox to make sure that the merged network state is feasible and that there is no modeling inconsistencies between our OPF and standard load-flow.

Three European situations are used:
\begin{itemize}
\item \textbf{09h30}: Data from the 12$^{th}$ February 2014 at 09h30
\item \textbf{11h15}: Data from the 12$^{th}$ February 2014 at 11h15
\item \textbf{08h30}: Data from the 13$^{th}$ February 2014 at 08h30
\end{itemize}

Some of the snapshots data are missing and DACF data are used to replace them in order to have the full European network:

\begin{center}
\begin{tabular}{l|c|c c c}
Country & Code & 09h30 & 11h15 & 08h30\\
 & & status & status & status\\\hline
Austria & AT & SN & DACF & SN\\
Belgium & BE & SN & SN & SN\\
Switzerland & CH & SN & SN & SN\\
Czech republic & CZ & DACF & SN & SN\\
Germany & DE & Both & Both & Both\\
France & FR & DACF & SN & SN\\
Croatia & HR & DACF & DACF & DACF\\
Italy & IT & SN & SN & SN\\
Netherlands & NL & SN & SN & SN\\
Poland & PL & DACF & DACF & DACF\\
Slovenia & SI & SN & SN & DACF
\end{tabular}
\end{center}
German network data are snapshot on some parts of the country and DACF on the others. Furthermore, due to TSO actions on the network, the grid topology may change over time, from situation to another.

\vskip.3cm
The indicators presented in the table below show the interest of the hierarchical merging procedure. This method manages to hold targets for active power flow on interconnections and voltage levels of PV buses. Deviations are mainly focused on loads which corresponds to the order of objectives in the hierarchical procedure.  

\begin{center}
\begin{tabular}{l|l|l|l}

Indicator & 09h30 & 11h15 & 08h30 \\\hline
\textit{Interco SN\_SN}   & 3  $/$ 68  &3  $/$ 99   &3  $/$ 96 \\
\textit{Interco SN\_DACF} & 24 $/$ 73  &14 $/$ 47   &4 $/$ 48   \\
\textit{V SN}             & 11 $/$ 242 &20 $/$ 455  &24 $/$ 451 \\
\textit{V DACF}           & 19 $/$ 329  &16 $/$ 73   &9 $/$ 82   \\
\textit{Loads ACT SN}     & 26 $/$ 1914&43 $/$ 2859 &39 $/$ 2806 \\
\textit{Loads REA SN}     & 14 $/$ 1884&40 $/$ 2848 &34 $/$ 2797 \\
\end{tabular}
\end{center}

\begin{figure}[!h]
\centering
\includegraphics[width=3.5in]{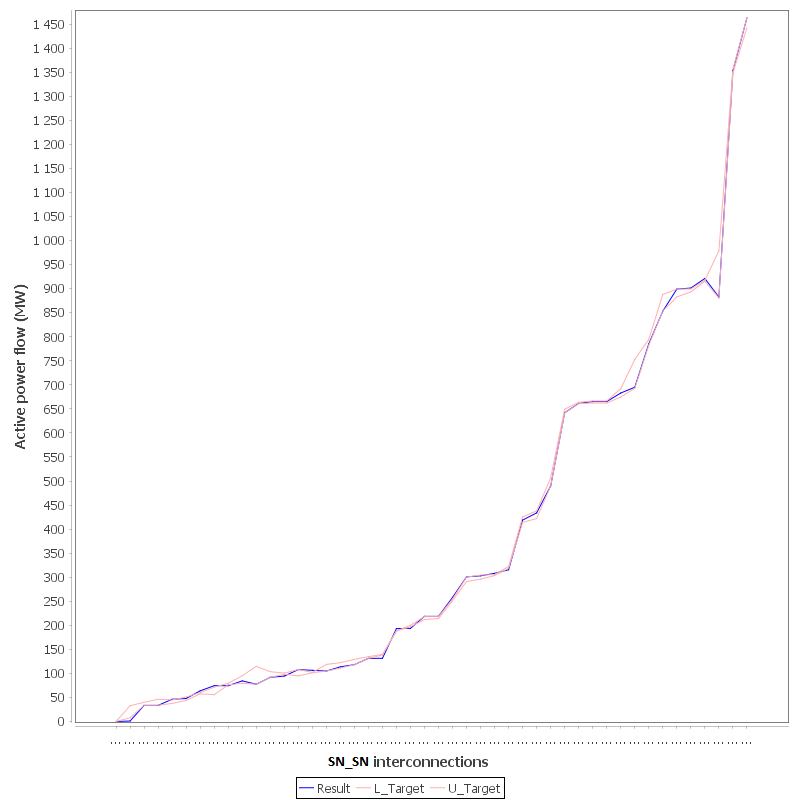}
\caption{Deviations of active power flow on SN\_SN interconnections for the 09h30 situation}
\label{09h30_SN_SN}
\end{figure}

\begin{figure}[!h]
\centering
\includegraphics[width=3.5in]{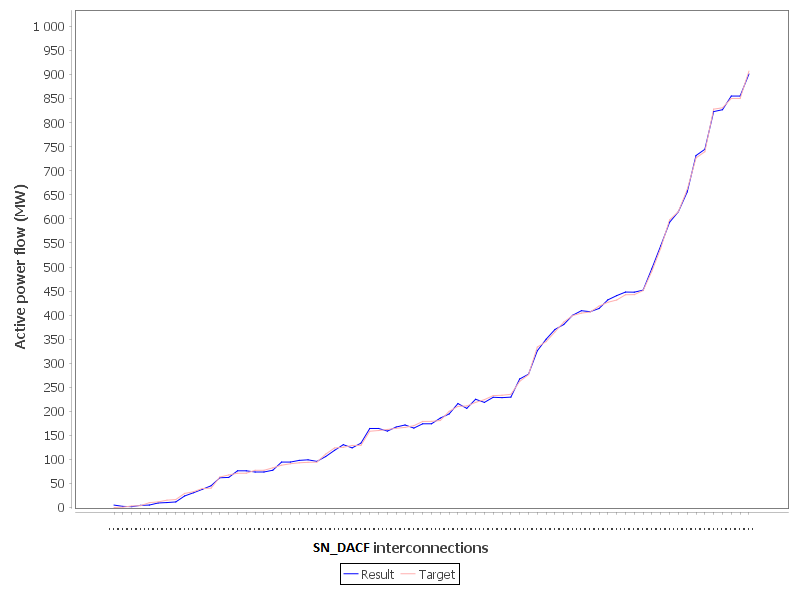}
\caption{Deviations of active power flow on SN\_DACF interconnections for the 09h30 situation}
\label{09h30_SN_DACF}
\end{figure}

\begin{figure}[!h]
\centering
\includegraphics[width=3.5in]{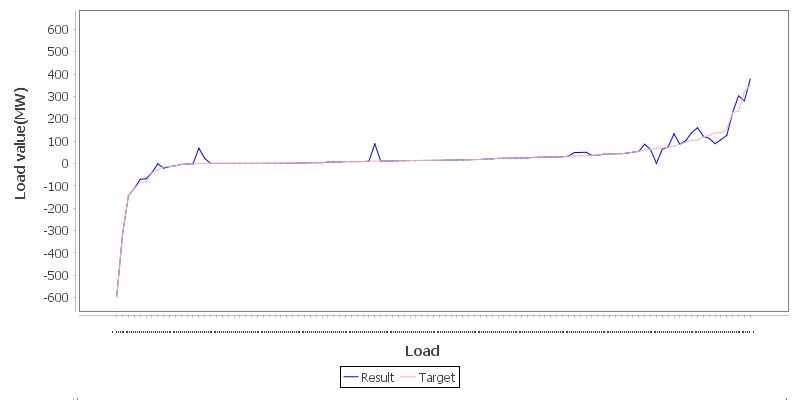}
\caption{Deviations of active SN loads for the 09h30 situation}
\label{09h30_LOAD_SN}
\end{figure}

For the 09h30 situation, 3 SN\_SN interconnections over 68  have deviations greater than $5MW$. Although, these deviations remains small. 
In fact, the maximum deviation from the two targets of SN\_SN interconnections (L\_TARGET and U\_TARGET) is $7MW$ and $9MW$ from SN\_DACF interconnections with a target from the SN side. 
These deviations are relatively small compared to active power flow on interconnection which reaches $1450MW$ as we can see in figures~\ref{09h30_SN_SN} and~\ref{09h30_SN_DACF}.
These two figures represent active power flows on interconnections ordered in function of target values. One can see that the active power flows remains very close to the target values. Concerning the voltage level, the maximum deviation is $0.37kV$ located on a DACF bus. Only 11 SN PV buses have a deviation larger than $0.01kV$. Regarding loads, the deviations of active load for SN buses (figure~\ref{09h30_LOAD_SN}) are up to $77MW$. Even if the maximum deviation is high, this is the expected behavior when using the $L_1$ norm. Indeed, deviations are local and have a large magnitude.

\vskip.3cm
For the 08h30 and 11h15 situations, the analysis is similar to 09h30 situation. The indicators are proportionally comparable to the 09h30 situation when taking into account the grid topology.
In the 11h15 situation, the maximum deviation of active power flow on interconnections is equal to $6MW$ for SN\_SN interconnections and $12 MW$ for SN\_DACF interconnections. Small deviations are detected on voltage level where the maximum value is $0.14kV$.
In the 08h30 situation, maximum deviation of active power flow on interconnections is equal to $8MW$ for SN\_SN interconnections and $6 MW$ for SN\_DACF interconnections while deviations of voltage level are up to $0.25kV$.

\subsection{Robustness against data errors}
Second experiment intends to test the robustness of the hierarchical approach against data errors. The aim of this experiment is to make sure that the introduction of DACF data, poorly reliable, does not affect the merging of the SN part of the network.

The 09h30 situation presented in~\ref{subsec:first_exp} is used for this experiment. Two situations are compared: the initial situation ($09h30$) and the situation where SN data from Netherlands (NL) is replaced by DACF data ($09h30$ $modified$). Hence, 3 parts of the network are distinguished: the $SN$ part in the two situations, the $DACF$ part in the two situations and the $NL$ network which SN in $09h30$ and DACF in $09h30$ $modified$. These 3 parts of the network are used for computing our indicators.

The indicators presented in the following table show that SN part of the network is not affected by the introduction of DACF data.

\begin{center}
\begin{tabular}{l|l|l}

Indicator & 09h30 & 09h30 modified \\\hline
\textit{V SN}             & 11 $/$ 231 & 44 $/$ 231 \\
\textit{V NL}             & 0 $/$ 11  & 0 $/$ 20 \\
\textit{Loads ACT SN}     & 23 $/$ 1284 & 22 $/$ 1284 \\
\textit{Loads ACT NL}     & 3 $/$ 630 & 6 $/$ 273 \\
\textit{Loads REA SN}     & 12 $/$ 1229& 11  $/$ 1229 \\
\textit{Loads REA NL}   & 2 $/$ 655& 13 $/$ 275 \\
\end{tabular}
\end{center}

Interconnections are not impacted by the introduction of DACF data. The number of deviations of active power flow on interconnections besides NL interconnections is the same for the two situations. The number of deviations of voltage levels for SN buses is more important in the modified situation but maximum deviation remains comparable with the initial situation ($0.13kV$ in $09h30$ and $0.15kV$ in $09h30$ $modified$).
Regarding loads, results show that the number of deviations of active load and reactive load for SN buses remains stable between the two situations. We also note that the network topology of the NL DACF data is different from the one used in the initial situation (e.g. 11 PV buses in $09h30$ and 20 in $09h30$ $modified$). Moreover, when NL data are DACF the number of deviations in the NL part of the grid is more important proportionally to the size of the NL network.

This experiment shows that the hierarchical merging procedure successfully consolidates the data even when some snapshot data are missing and replaced with DACF data. The use of the DACF data, less reliable than the SN data does not impact the consolidation in the SN part of the network.

\subsection{Numerical robustness}
Third experiment aims to check the robustness of the hierarchical approach compared to the merging procedure used at CORESO.
The merging tool deployed at CORESO usually fails less than 1\% of the time due to numerical difficulties. On a week of July 2015, the number of failures have reached 5\%. The set of situations used in this third experiment corresponds to this week of data: 672 situations which corresponds to a situation, every 15 minutes of the 7 days.
Results of this experiment shows that the hierarchical merging approach is as reliable as the merging procedure used at CORESO.

\section{Conclusion}
In this paper, a hierarchical merging procedure solving several successive ACOPFs is introduced to build a consistent merged network state by minimizing the deviations from the data provided by the different TSOs. The aim of this approach is to make the right trade-off between the different objectives defined in function of deviation type (voltage levels, active power flows or loads) and data reliability level (SN or DACF). The experiments on the European network have shown the interest of this approach to obtain a consistent feasible merged network state. They also demonstrate the robustness of such approach when dealing with data errors. The quality of the DACF data (used as replacement when SN data are not available) does not influence the merging quality on the SN part of the network. Moreover, the method have been compared to the merging procedure used in CORESO. The results exhibit that the hierarchical merging procedure is reliable.

Further work may improve the hierarchical merging procedure by detecting Phase-Shifting Transformer configurations errors and correct them using an optimization process. Another way of improvement could be the use of the previous merged network state as a replacement of SN data or as an initial point from the optimization process.

\bibliographystyle{IEEEtran}
\bibliography{biblio}

\end{document}